\documentclass{amsart}
\usepackage{amsfonts,amsmath,amsthm,amssymb}

\newtheorem{lemma}{Lemma}
\newtheorem{theorem}{Theorem}
\newtheorem{corollary}{Corollary}

\def\gp#1{\langle #1 \rangle}
\def\m1{^{-1}}

\title{Lie properties of crossed products}
\author{Adalbert Bovdi and Alexander Grishkov}

\address{
A.~Bovdi\newline Institute of Mathematics, University of
Debrecen,\newline P.O.  Box 12, H-4010 Debrecen, Hungary}
\email{bodibela@math.klte.hu}

\address{
A.~Grishkov
\newline Departamento de Matem\'atica, (IME-USP),
\newline  Rua do Matao, 1010 - Cidade Universit\'aria,
\newline  CEP 05508-090, Sao Paulo - SP,  Brasil}
\email{grishkov@ime.usb.br}

\dedicatory{Dedicated to 60th birthday of Professor
I.P.~Shestakov}

\thanks{Supported by OTKA  No.K68383, by FAPESP (proc. 06-56203-3) and by RFFI(07-01-00392A)}

\subjclass {Primary: 16S34, 16U60; Secondary: 20C05}

\keywords{twisted group algebra, group of units, crossed product,
Lie nilpotent}

\begin{document}

\begin{abstract}
Let $F^\lambda_{\sigma} [G]$ be a crossed product of a group $G$
and the field $F$. We study the Lie properties of
$F^\lambda_{\sigma} [G]$ in order to obtain a characterization of
those crossed products which are upper (lower) Lie nilpotent and
Lie $(n,m)$-Engel.
\end{abstract}

\maketitle

\section{Introduction}

Let $G$ be a group and let $Aut(F)$ be the group of automorphisms
of a field $F$. Assume the mapping $ \sigma: G\to Aut(F)$ and
the twisting function  $ \lambda: G\times G\to U(F) $
satisfy the following conditions:
\begin{equation}\label{E:1}
\lambda(a, bc)\lambda (b, c)=\lambda (ab, c)\lambda (
a,b)^{\sigma(c)}
\end{equation} and
\begin{equation}\label{E:2}
\alpha ^{\sigma(a)\cdot \sigma(b)}=  \alpha^{\sigma(ab)}
\end{equation}
for all $a, b, c\in G$ and $\alpha(a,b)\in F$. The twisting
function $\lambda$ is also called  a {\bf factor system} of the
group $G$ over the field $F$ relative to $\sigma$. It is  a
$2$-cocycle of the group  of units $U(F)$ of $F$ with the natural
$G$-module structure (in the cohomology group of  $G$ in $U(F)$).

Assign to every  $g\in G$  a symbol $\widetilde{g}$, and let
$$\textstyle
F^\lambda_{\sigma}[G]=\{\sum_{g\in G}\widetilde{g}\alpha_g \mid
\alpha_g\in F\}
$$
be the set of all formal sums with  finitely many nonzero
coefficients $\alpha_g$. Two elements $\textstyle x=\sum_{g\in G}
\widetilde{g}\alpha_g \quad\text{and} \quad y=\sum_{g\in
G}\widetilde{g}\beta_g $ from $F^\lambda_{\sigma} [G]$ are equal
if and only if $\alpha_g =\beta_g$ for all $g\in G$. On the set
$F^\lambda_{\sigma} [G]$ addition and multiplication are defined
as follows:
\begin{itemize}
\item[(i)] $\sum_{g\in G}\widetilde{g}\alpha_g + \sum_{g\in
G}\widetilde{ g}\beta_g=\sum_{g\in G}
\widetilde{g}(\alpha_g+\beta_g);$ \item[(ii)] $\widetilde{
g}\widetilde{h}=\widetilde{gh}\lambda (g, h)$, where  $\lambda$ is
the twisting function; \item[(iii)] $\alpha \widetilde{g}=
\widetilde{g} \alpha^{\sigma(g)}$.
\end{itemize}
The product of arbitrary elements $x$ and $y$ is  determined by
distri\-bu\-tivity. It is easy to check that $F^\lambda_{\sigma}
[G]$ is an associative ring which  is called a {\bf crossed
product} of the group $G$ over the field $F$. If $\sigma$ is a
trivial mapping then we shall denote this  ring by $F^\lambda [G]$
and it is a {\bf twisted group algebra}.

From (\ref{E:1}) and  (\ref{E:2}) follows that
$$
\lambda (h, 1)^{\sigma(1)}= \lambda (h, 1)= \lambda (1, 1), \qquad
\lambda (1, h)=\lambda (1, 1)^{\sigma(h)}.
$$
Clearly,\quad  $\widetilde{1}\cdot \lambda (1, 1)^{-1}$\quad  is
the identity element of the ring $F^\lambda_{\sigma} [G]$ and
\begin{equation}\label{E:3}
\widetilde{g}^{-1}=\lambda (g^{-1}, g)^{-1}\lambda (1,
1)^{-1}\widetilde{g^{-1}}= \widetilde{g^{-1}}\lambda (g,
g^{-1})^{-1}\lambda (1, 1)^{-1}.
\end{equation}
For the twisted group algebra $F^\lambda [G]$ we can always assume
without loss of generality, that the twisting  function is
normalized, that is
$$
\lambda (h, 1)= \lambda (1, 1)= \lambda (1, h)=1,
$$
so $\widetilde{1}$ is its identity element.

In the  theory of ordinary group rings the Lie properties play an
important role. Group algebras with the many "good"  Lie
properties were   described during the 70's using the theory on
polynomial identities. Later these results  were applied to the
study of the group of units. In  most cases the Lie structure
reflects very  well the characteristics of the group of units,
there is a close relationship between their properties of them.
Here we  describe   the structure of those crossed products  which
are upper (lower) Lie nilpotent and Lie $(n,m)$-Engel. We
generalize results of Passi, Passman and Sehgal \cite{pps} which
were obtained for the group algebras.

\section{Preliminaries Results}

Let $F^\lambda [G]$ be a twisted group algebra  with normalized
twisting  function $\lambda$. Then
$$
W=\{w\in G\mid \lambda(g,w)=\lambda(w,g)=1 \;\;\text{for all}\;\;
g\in G\}
$$
is a subgroup of $G$. Indeed, if $w_1,w_2, w\in W$ and $g\in G$,
then
\[
\begin{split}
\widetilde{(w_1w_2)}\widetilde{g}&=\widetilde{w_1}\widetilde{w_2}\widetilde{g}=\widetilde{w_1}(\widetilde{w_2g})=\widetilde{w_1w_2g};\\
\widetilde{g}\widetilde{(w_1w_2)}&=\widetilde{g}\widetilde{w_1}\widetilde{w_2}=
\widetilde{(gw_1)}\widetilde{w_2}=\widetilde{gw_1w_2}.
\end{split}
\]
So $\lambda(w_1w_2,g)=\lambda(g,w_1w_2)=1$ and $w_1w_2\in W.$
Furthermore, by (\ref{E:1}) and (\ref{E:3}), it is easy to check that
\[
\begin{split}
\widetilde{w^{-1}}\widetilde{g}=\widetilde{w}^{-1}\widetilde{g}
=(\widetilde{g}^{-1}\widetilde{w})^{-1}&=\lambda(g,g^{-1})(\widetilde{g^{-1}}\widetilde{w})^{-1}
=\lambda(g,g^{-1})(\widetilde{g^{-1}w})^{-1}\\
&=\lambda(g,g^{-1})\lambda(g^{-1}w,w^{-1}g)^{-1}\widetilde{w^{-1}g}\\
&=\lambda(g^{-1},w)^{-1}\lambda(w,w^{-1}g)^{-1}\widetilde{w^{-1}g}
= \widetilde{w^{-1}g}.
\end{split}
\]
Similarly, $\widetilde{g}\widetilde{w^{-1}}= \widetilde{gw^{-1}}$,
and this shows that  $w^{-1}\in W$, so $W$ is a subgroup.

Let $H$ be  a normal subgroup of $G$ such that
$H\subseteq W$ and denote by  ${\bf I}(H)$  the ideal of
$F^\lambda[G]$ generated by the elements
$\widetilde{h}-\widetilde{1}$ with $h \in H$. It is easy to see
that if $\{u_i\}$ is a transversal of $H$ in $G$ then the elements
of the form\;  $\widetilde{u_i}(\widetilde{h}-\widetilde{1})$ \;
with $h\ne 1$ constitute an  $F$-basis of the ideal ${\bf I}(H)$.

Define a new  function: $\mu: G/H \times G/H \to U(F)$ the
following way:  Choose two representatives  $g_i=u_kh_1$ and
$g_j=u_lh_2$ with $h_i\in H\subseteq W$ of  the cosets of $H$ in
$G$.   Then by (\ref{E:1}) in  $F^\lambda[G]$ we get
\[
\begin{split}
\lambda(g_i,g_j)&=\lambda(u_kh_1,u_lh_2)\\
&=\lambda(u_k,h_1)^{-1}\lambda(h_1,u_lh_2)\lambda(u_k,h_1u_lh_2)=\lambda(u_k,h_1u_lh_2)
\end{split}
\]
and  $h_1u_lh_2=u_lh$ for suitable $h\in H$. It is easy to check
that
\[
\begin{split}
\lambda(u_k,h_1u_lh_2)&= \lambda(u_k,u_lh)\\
&=\lambda(u_l,h)^{-1}\lambda(u_k,u_l)\lambda(u_ku_l,h)=\lambda(u_k,u_l).
\end{split}
\]
We proved that the twisting function $\lambda$ satisfies the
condition $\lambda(g_i,g_j)=\lambda(u_k,u_l)$ and it can  define
the function
\begin{equation}\label{E:4}
\mu(g_iH,g_jH)=\lambda(g_i,g_j),
\end{equation}
which is a twisting  function of $G/H$.  Of course, any element of
$F^\lambda [G]$ has the form $\sum_iu_ix_i$ with $x_i\in F^\lambda
[H]$ and $x_i+ {\bf I}(H)=\lambda_i + {\bf I}(H)$ for suitable
$\lambda_i\in F$. Now it is  easy to see that $F^\lambda [G]/{\bf
I}(H)$ is an $F$-algebra with a basis consisting  of the images
$\widetilde{u_i}$ of the  coset representatives $u_i$ of $H$ in
$G$. We obtain the isomorphism:
\begin{equation}\label{E:24}
F^\lambda [G]/{\bf I}(H)\cong F^\mu [G/H].
\end{equation}
Unlike  ordinary group algebras, the twisted group algebra
$F^\lambda [G]$ does not have a natural group basis. If
the algebra $F^\lambda [G]$ has an $F$-basis
$\widetilde{G}=\{\widetilde{g}\mid g\in G\}$ such that for each
$g$ of $G$ there exists an  element $d_g\in F$ such that the
elements of the set $\{d_g\widetilde{g}\mid g\in G\}$ form a group
basis for $F^\lambda [G]$, then  $F^\lambda [G]$ is called  {\bf
untwisted}. In this situation  $F^\lambda [G]$ is isomorphic to
$FG$ via this diagonal change of the basis. In addition,
$F^\lambda [G]$ is called {\bf stably untwisted} if  there exists
an extension $K$ of the field $ F$ such that $K^\lambda [G]$ is
untwisted.

We use a criteria from \cite{passm}  to verify that a twisted
group algebra is  untwisted or  stably untwisted.
\begin{lemma}\cite{passm}\label{L:1}
Let $F^\lambda [G]$ be a twisted group algebra.
\begin{itemize}
\item[(i)] $F^\lambda [G]$ is  untwisted if and only if there
exists an $F$-algebra homomorphism  $F^\lambda [G]\to F$.
\item[(ii)]  $F^\lambda [G]$ is  stably untwisted if and only if
there  exists an $F$-algebra homomorphism of $F^\lambda [G]$ into
a commutative $F$-algebra $R$.
\end{itemize}
\end{lemma}

\bigskip
Recall that  $\{u_i\}$ is a transversal of $H$ in $G$ and  each
element from $G$ can be written uniquely  in  the form $g=zu_i$
with $z\in H$.

\begin{lemma}\label{L:5} Let $F^\lambda [G]$ be a twisted group
algebra and assume that the twisting  function $\lambda$ of the
normal subgroup $H$ of $G$ satisfies the condition
$\lambda(h_1,h_2)=1$ for all $h_1,h_2\in H$. Then the algebra
$F^\lambda [G]$ can be realized alternatively  as a twisted group
algebra $F^\tau [G]$ with the following diagonal change of the
basis by
$$
\overline{g}= \begin{cases}
\widetilde{g},                 & \quad  \text{if }\quad  g\in H;\\
\lambda(z,u_j)\widetilde{g} & \quad  \text{if }\quad  g=zu_j\in G\setminus
H,\quad z\in H,
\end{cases}
$$
with twisting function  $\tau$. This twisting  function
has the property $\tau(h,g)=1$ for all $h\in H$ and $g\in G$.

If  $\tau(g,h)=1$ for all $h\in H$ and $g\in G$,  then
$$
F^\tau [G]/{\bf I}(H)\cong F^\tau [G/H].
$$
\end{lemma}
 {\bf Proof}. It is easy to see that if $g=h_1u_j$ then
\[
\begin{split}
\tau(h,g)\overline{hg}&=\overline{h}\cdot
\overline{h_1u_j}=\widetilde{h}\cdot
\widetilde{h_1u_j}\lambda(h_1,u_j)\\
&= \widetilde{h} \cdot\widetilde{h_1}\cdot\widetilde{u_j}=
\widetilde{hh_1}\cdot\widetilde{u_j}=
\lambda(hh_1,u_j)\widetilde{hg}=\overline{hg}.
\end{split}
\]
This yields \quad   $\tau(h,g)=1$.

If $\tau(g,h)=1$ for all  $h\in H$ and  $g\in G$,  then we can
apply (\ref{E:24}) to obtain the lemma.
\bigskip

The twisted group algebras satisfying polynomial identities  have
been classified in \cite{pass} and this result was modified in
\cite{liu}  for the  stably untwisted case.
\begin{lemma}\cite{liu} \label{L:11}
 If $F^\lambda[G]$ is a twisted group algebra of positive
characteristic $p$ with  a polynomial identity of degree $n$ then
$G$ has a subgroup $A$ of finite index  such that $F^\lambda[A]$
is  stably untwisted, the commutator subgroup $A'$ of $A$ is a
finite $p$-group and $|G:A||A'|$ is bounded by a fixed function of
$n$.
\end{lemma}

\smallskip
Let  $F^\lambda[G]$ be a twisted group algebra. For each $h\in G$
of order $k$ we have
\begin{equation}\label{E:8}
\textstyle\widetilde{h}^{k}=\prod_{i=1}^{k-1}\lambda(h^i,h)\cdot
\widetilde{1}.
\end{equation}
It is convenient  to say that \quad
$\mu(h)=\prod_{i=1}^{k-1}\lambda(h^i,h)$\quad is the {\bf twist}
of $\widetilde{h}$ which plays  an important role in the study of
$F^\lambda[G]$. An   $p$-element $u$ of $G$ is called an {\bf
untwisted $p$-element} if the order of $u$ coincides with the
order of $\widetilde{u}\gamma$ for some $\gamma\in F$.

\begin{lemma}\label{L:2}
Let  $F^\lambda [G]$ be a twisted group algebra of positive
characteristic $p$ such  that the commutator ideal $F^\lambda
[G]^{[2]}=[F^\lambda [G],F^\lambda [G]]F^\lambda [G]$ is a nil
ideal and let $a,b\in G$.
\begin{itemize}
\item[(i)] If $\widetilde{a}$ is a $p$-element, then $a$ is also a
$p$-element, its order coincides with the order of $\widetilde{a}$
and its  twist is $\mu(a)=1.$ \item[(ii)] If $a$ and $b$ are
untwisted $p$-elements, and the orders of $\widetilde{a}\gamma_1$
and $\widetilde{b}\gamma_2$ coincide with the order of $a$ and $b$
respectively for some $\gamma_i\in F$  then $ab$ is also a
$p$-element, the order $p^l$ of $ab$ coincides with the
order of $\widetilde{a}\widetilde{b}\gamma_1\gamma_2$
and
$$
\mu(ab)=(\gamma_1\gamma_2\lambda(a,b))^{-p^l}.
$$
\item[(iii)] The group commutator $(a,b)$  is an untwisted
$p$-element for all $a,b\in G$ and if $p^m$ is the order of
$(a,b)$, then
$$ \bigl(\lambda(a,b)^{-1}\lambda(b,b^{-1}ab)\lambda(a,(a,b)
\lambda((b,a),(a,b))^{-1}\bigr)^{-p^m}=\mu((a,b)).
$$
\end{itemize}
\end{lemma}

{\bf Proof}.  $(i)$\quad Let  $\widetilde{a}$ be an element of
order $p^t$ and let $k$ be the order of $a$. Then  (\ref{E:8})
yields that $\widetilde{a}^{k}=\mu(a)\widetilde{1}$, so
$\widetilde{a}^{kp^t}=\mu(a)^{p^t}\widetilde{1}=\widetilde{1}$. It
follows that the twist $\mu(a)$ of $\widetilde{a}$ satisfies the
condition  $\mu(a)^{p^{t}}=1$.  But in the field $F$ of
characteristic $p$ this is possible only if $\mu(a)=1$. Hence the
order of $a$ coincides with the order of $\widetilde{a}$ and
$\mu(a)=1$.

$(ii)$\quad First note that if $F^\lambda [G]^{[2]}$ is a nil
ideal, then
$$
xy\equiv yx \pmod{F^\lambda [G]^{[2]}},\qquad (xy)^{n}\equiv x^ny^n
\pmod{F^\lambda [G]^{[2]}}
$$
for all $x,y\in F^\lambda [G]$ and for all $n$. It follows that
the nilpotent elements of $F^\lambda [G]$ form an ideal $N$ and
$F^\lambda [G]/N$ is commutative.

Let $a$ and $b$ be untwisted $p$-elements, and assume that the
orders of $\widetilde{a}\gamma_1$ and $\widetilde{b}\gamma_2$
coincide with the order of $a$ and $b$  respectively for some
$\gamma_i\in F.$ Since
$\widetilde{a}\gamma_1-\widetilde{1}\in N$ and
$\widetilde{b}\gamma_2-\widetilde{1}\in N$, we have
$$
\widetilde{a}\widetilde{b}\gamma_1\gamma_2-\widetilde{1}=(\widetilde{a}\gamma_1-\widetilde{1})
(\widetilde{b}\gamma_2-
\widetilde{1})+(\widetilde{a}\gamma_1-\widetilde{1})+(\widetilde{b}\gamma_2-\widetilde{1})\in
N,
$$
because the  nilpotent elements of $F^\lambda [G]$ form an ideal.
It follows that
$$
(\widetilde{a}\widetilde{b}\gamma_1\gamma_2-\widetilde{1})^{p^{l_1}}=0
$$
for some $l_1$, which shows that
$\widetilde{a}\widetilde{b}\gamma_1\gamma_2$ has order $p^l$. Now
$(\widetilde{ab}\lambda(a,b)\gamma_1\gamma_2)^{p^l}=\widetilde{1}
$ implies that $ab$ has order $p^m$ and,  by  (\ref{E:8}), we
obtain that
$$
(\widetilde{ab}\lambda(a,b)\gamma_1\gamma_2)^{p^m}=
\mu(ab)\cdot(\lambda(a,b)\gamma_1\gamma_2)^{p^m}\cdot\widetilde{1}
$$
which yields\;
$(\mu(ab)\cdot(\lambda(a,b)\gamma_1\gamma_2)^{p^m})^{p^{l-m}}=1$.\;
Then $m=l$ and
$$
\mu(ab)(\lambda(a,b)\gamma_1\gamma_2)^{p^m}=1.
$$
Hence $ab$ has
order $p^l$ and $\mu(ab)= (\lambda(a,b)\gamma_1\gamma_2)^{-p^m}$.

$(iii)$ Obviously,  the Lie commutator
$[\widetilde{a},\widetilde{b}]$ belongs to the nil ideal
$F^\lambda [G]^{[2]}$ for all $a,b\in G$ and the identity
\begin{equation}\label{E:9}
[\widetilde{a},\widetilde{b}]=\widetilde{a}^{-1}\widetilde{b}^{-1}
((\widetilde{a},\widetilde{b})-\widetilde{1})
\end{equation}
ensures   that  $(\widetilde{a},\widetilde{b})-\widetilde{1}$ is
nilpotent. Then  $(\widetilde{a},\widetilde{b})$ has order $p^m$
and an easy  computation shows  that
\begin{equation}\label{E:10}
(\widetilde{a},\widetilde{b})=\widetilde{(a,b)}\chi((a,b))
\end{equation}
for some $ \chi((a,b))$ of $F$. The argument of the proof of (i)
states that $p^m$ is the order of $(a,b)$ and
\begin{equation}\label{E:11}
(\widetilde{a},\widetilde{b})^{p^m}=\mu((a,b))
\chi((a,b))^{p^m}\widetilde{1}=\widetilde{1}.
\end{equation}
 Moreover,  from (\ref{E:3})
and  (\ref{E:10})  follows that
\[
\begin{split}
\chi((a,b)) =&\widetilde{a}^{-1}\widetilde{b}^{-1}
\widetilde{a}\cdot\widetilde{a^{-1}ba}\lambda(b, (b,a))
\lambda((b,a),(a,b))^{-1}\\
=&\widetilde{a}^{-1}\widetilde{b}^{-1}\widetilde{ba}
\lambda(a,a^{-1}ba)\lambda(b,(b,a)) \lambda((b,a),(a,b)^{-1}\\
=&\lambda(b,a)^{-1}\lambda(a,a^{-1}ba)\lambda(b,(b,a))
\lambda((b,a),(a,b))^{-1}.
\end{split}
\]

\begin{corollary}\label{C:3} If  $F^\lambda[G]$ is  a twisted group algebra of
positive characteristic $p$ such  that $F^\lambda [G]^{[2]}$ is a
nil ideal then  the group commutator $(a,b)$ is an untwisted
$p$-element for any $a,b\in G$ and the product of  untwisted
$p$-elements of $G$ is also an untwisted  $p$-element.
\end{corollary}

\bigskip

\section{ Upper and lower Lie nilpotent\\ crossed products}

Let $R$ be an associative ring. The lower Lie central  series in
$R$ is defined   inductively as follows:
\[
\gamma_1(R)=R, \quad  \gamma_2(R)=[\gamma_1(R),  R], \ldots,
\gamma_n(R)=[\gamma_{n-1}(R), R], \ldots .
\]
The two-sided ideal $ R^{[n]}=\gamma_n(R)R$ of $R$ is called the
{\bf  $n$-th lower Lie power} of  $R$. By Gupta and Levin
\cite{gupt}, these ideals satisfy the conditions
\begin{equation}\label{E:5}
R^{[m]}R^{[n]}\subseteq R^{[n+m-2]},\qquad (n, m\geq 2).
\end{equation}

Let us  define by induction a second set of ideals  in  $R$:
$$
R^{(1)}=R, \quad  R^{(2)}=[R^{(1)}, R]R,   \ldots,
R^{(n)}=[R^{(n-1)}, R]R, \ldots
$$
The ideal $R^{(n)} $ is called {\bf the $n$-th upper Lie power} of
$R$ and these ideals have the   property
\begin{equation}\label{E:6}
R^{(n)}R^{(m)}\subseteq R^{(n+m-1)}, \qquad (n, m\geq 1).
\end{equation}

Recall that $R$ is called   {\bf upper Lie nilpotent} if
$R^{(n)}=0$ for some $n$. Similarly, the ring $R$ with $R^{[m]}=0$
for some $m$ is called {\bf lower Lie nilpotent}; these classes of
rings are different.

First assume that $R$ is lower Lie nilpotent and let $R^{[t]}=0$.
For $k\geq 3$ we  choose an $l$ such that $l(k-1)+2\geq t$. Then
(\ref{E:5}) forces  $ (R^{[k]})^l= R^{[l(k-1)+2]}=0, $ so
$R^{[k]}$ \quad $(k\geq 3)$ is a nilpotent ideal. Note that every
element of $R^{[2]}$ has the form
$$
x=[a_1,b_1]r_1+[a_2,b_2]r_2+\cdots+[a_s,b_s]r_s
$$
for some $a_i,b_i,r_i\in R$ and
\begin{equation}\label{E:7}
\begin{split}
([a_j,b_j]r_j)^2=&r_j[a_jb_j,b_j,a_j]r_j+[a_jb_j,b_j,a_j]r_j\\
&+[a_j,b_j,a_j]b_jr_j+ [a_jb_j,b_j,r_j][a_jb_j,b_j]r_j.
\end{split}
\end{equation}

Let $R$ be of characteristic $p$. Then, by Brauer's formula,
$$
x^p=([a_1,b_1]r_1)^p+([a_2,b_2]r_2)^p+\cdots+([a_s,b_s]r_s)^p+z
$$
for suitable  $z=[c_1,d_1]+[c_2,d_2]+\cdots+[c_q,d_q]\in [R,R]$
and $c_i,d_i\in R$. Now  (\ref{E:7})   ensures that
$([a_j,b_j]r_j)^p\in R^{[3]}$, because $R^{[3]}$ is an ideal.
Since the elements $[c_i,d_i]$ and $[c_j,d_j]$ commute  modulo
$R^{[3]}$,  Brauer's formula implies that
$$
x^{p^2}=([c_1,d_1])^p+([c_2,d_2])^p+\cdots+([c_q,d_q])^p\in
R^{[3]}.
$$
There exists $s$ such that $x^{p^s}=0$ for every $x\in R^{[2]}$,
because $R^{[3]}$ is a nilpotent ideal, so the ideal $R^{[2]}$ is
nil.

Similar results are  valid for upper Lie nilpotent rings  $R$.

\begin{lemma}\label{L:3}
Let  $F^\lambda [G]$ be a twisted group algebra of positive
characteristic $p$ such that either\quad  $F^\lambda [G]^{[m]}=0$,
or $F^\lambda [G]^{(m)}=0$.\quad  If $p^t\geq m$,  then
\begin{itemize}
\item[(i)] $b^{p^t}$ is a central element  of $G$ for any $b\in
G$; \item[(ii)] for  $q\ne p$ each $p$-element $a\in G$ commutes
with any $q$-element $c\in G$ and
$$
\lambda(a,c)=\lambda(c,a).
$$
\end{itemize}
\end{lemma}

{\bf Proof.}\quad For  $a$ and $b$  of $G$ the well known Lie
commutator formula confirms
 $$
[\widetilde{a},\widetilde{b},p^t]=[\widetilde{a},
\underbrace{\widetilde{b},\widetilde{b},\ldots,
\widetilde{b}}_{p^t}]=[\widetilde{a},\widetilde{b}^{p^t}]=0.
$$
Hence
$\widetilde{a}\widetilde{b}^{p^t}=\widetilde{b}^{p^t}\widetilde{a}$,
this  yields    $ab^{p^t}=b^{p^t}a$ and
$\lambda(a,b^{p^t})=\lambda(b^{p^t},a)$ for all $a\in G$. Since
any $q$-element $c$ of $G$ can be written as $c=b^{p^t}$ for some
$b$, the desired assertion follows.

\bigskip
 We present  two different proofs for the next result.

\begin{lemma}\label{L:7} If $F^\lambda [G]$ is a twisted group algebra of $char(F)=p$
such that either $F^\lambda [G]^{[p^t]}=0$ or $F^\lambda
[G]^{(p^t)}=0$, then  $G'$ is a finite $p$-group.
\end{lemma}

{\bf Proof}. The first proof of the lemma uses  the theory of
polynomial identities.  As we showed before  $F^\lambda[G]$
satisfies  a polynomial identity, and by Lemma \ref{L:11}, the
group $G$ has a normal subgroup $A$ of finite index such that
$F^\lambda[A]$ is stably untwisted and  the commutator subgroup
$A'$ of $A$ is a finite $p$-group.

We start with some facts which will be used  freely.

1. If $g,h\in G$ and $(g,h)=1$, then
$[\widetilde{g},\widetilde{h}]=0$. Indeed, $
[\widetilde{g},\widetilde{h}]=\widetilde{gh}(\lambda(g,h)-\lambda(h,g))$
is nilpotent which is possible only if \quad
$\lambda(g,h)-\lambda(h,g)=0.$

2. If $h$ is  central in $G$, then $\widetilde{h}$ is a central
element of $F^\lambda [G]$.

3. By Lemma \ref{L:2},    $(a,b)$ is a $p$-element for all $a,b\in
G$ and

\begin{equation}\label{E:18}
(\widetilde{a},\widetilde{b})=\widetilde{(a,b)}\chi((a,b)).
\end{equation}
Now let $F^\lambda [G]$ be upper (lower) Lie nilpotent and assume
that $A$ is abelian. If $P$ is the maximal $p$-subgroup of $A$,
then $A=P\times Q$ for a suitable central $p'$-subgroup $Q$,
because $G'$ is $p$-group by Lemma \ref{L:2}. Moreover, $P$
belongs to the $FC$-center of $G$ and assume that
 $C_P(g)$ has infinite index in $P$ for some $g$ of $G$. For brevity, put
$\chi((g,g_i))=\pi_i$. Clearly, $(g,g_1)\ne 1$ for suitable
$g_1\in P$ and, using the fact that $F^\lambda[A]$ is commutative,
we have
$$
[\widetilde{g},\widetilde{g_1}]=\widetilde{g}\widetilde{g_1}(\widetilde{1}-(\widetilde{g_1},\widetilde{g}))
=\widetilde{g}\widetilde{g_1}(\widetilde{1}-\widetilde{(g_1,g)}\pi_1).
$$
Since  the subset $\{(h, g)\mid h\in P\}$ of $P$ is infinite,
there exists $g_2\in P$ such that
$$
(\widetilde{1}-\widetilde{(g_2,g)}\pi_2)(\widetilde{1}-\widetilde{(g_1,g)}\pi_1)\ne
0.
$$
Clearly,\quad  $[\widetilde{g_1},\widetilde{g_2}]=0$ and
\[
\begin{split}
[\widetilde{g},\widetilde{g_1},\widetilde{g_2}]=
[\widetilde{g}\widetilde{g_1}(\widetilde{1}-\widetilde{(g_1,g)}\pi_1),\widetilde{g_2}]
&=[\widetilde{g},\widetilde{g_2}]\widetilde{g_1}(\widetilde{1}-\widetilde{(g_1,g)}\pi_1)\\
&=\widetilde{g}\widetilde{g_2}\widetilde{g_1}(\widetilde{1}-\widetilde{(g_2,g)}
\pi_2)(\widetilde{1}-\widetilde{(g_1,g)}\pi_1).
\end{split}
\]
As before, it is easy to see that for each $n$ there exist
$g_1,g_2,\ldots, g_n$  in $P$ such that
$$
(\widetilde{1}-\widetilde{(g_n,g)}\pi_n)(\widetilde{1}-\widetilde{(g_{n-1},g)}\pi_{n-1})
\cdots(\widetilde{1}-\widetilde{(g_1,g)}\pi_1)\ne 0
$$
and
\[
\begin{split}
[\widetilde{g},&\widetilde{g_1},\ldots,\widetilde{g_n}]=\\
&=\widetilde{g}\widetilde{g_n}\widetilde{g_{n-1}}\cdots\widetilde{g_1}(\widetilde{1}-\widetilde{(g_n,g)}\pi_n)
(\widetilde{1}-\widetilde{(g_{n-1},g)}\pi_{n-1})\cdots(\widetilde{1}-\widetilde{(g_1,g)}\pi_1).
\end{split}
\]
Clearly $[\widetilde{g},\widetilde{g_1},\ldots,\widetilde{g_n}]\ne
0$ for each $n$, contradicting to the assumption that $F^\lambda
[G]$ is upper (lower) Lie nilpotent. Thus for any $g$ of $G$ the
centralizer $C_P(g)$ has finite index in $P$. But this imply that
$C_G(g)$ has finite index in $G$, because $Q$ is a central
subgroup and $A$ has finite index. Therefore we can suppose below
that $G$ is an $FC$-group.

Assume that the $p$-group $G'$ is infinite. Then $P$ is infinite,
$b_1=(g_1,g_2)\ne 1$ for suitable $g_1, g_2$ and
\[
\begin{split}
[\widetilde{g_1},\widetilde{g_2}]=\widetilde{g_2}\widetilde{g_1}((\widetilde{g_1},\widetilde{g_2})-\widetilde{1})
&=\widetilde{g_2}\widetilde{g_1}(\widetilde{(g_1,g_2)}-\widetilde{1})\\
&=\widetilde{g_2}\widetilde{g_1}(\chi((g_1,g_2))\widetilde{b_1}-\widetilde{1}).
\end{split}
\]
Now we claim that there exists a sequence $\{g_i\}$ with the
following properties: $ g_{2n+1}$, $g_{2n+2}$ and
$b_{2n+1}=(g_{2n+1}, g_{2n+2})$ are  elements of
$$
C_G(\{g_1,g_2,\ldots,g_{2n}, b_1,b_3,\ldots,
b_{2n-1} \})
$$
and
$$
(g_{2n+1},g_{2n+2})=b_{2n+1}\notin\langle b_1,b_3,\ldots,
b_{2n-1}\rangle.
$$
Indeed, assume that the sequence $g_1,g_2,\ldots,g_{2n}$  is
given. The subgroup
$$
C= C_G(\{g_1, g_2, \ldots, g_{2n},b_1,b_3,\ldots, b_{2n-1} \})
$$
of the $FC$-group $G$ has finite index.  Neumann's result
\cite{neu} asserts that in an $FC$-group with infinite commutator
subgroup $G'$   the commutator subgroup  of the   subgroup $C$ is
also infinite. Then $C'$ is a group with infinite commutator
subgroup and it contains elements $g_{2n+1}$ and $g_{2n+2}$ such
that
$$
b_{2n+1}=(g_{2n+1}, g_{2n+2}) \not\in \langle b_1, b_3, \ldots,
b_{2n-1} \rangle,
$$
because   $G'$ is a locally finite $p$-subgroup. Recall that if
$(a,b)=1$, then $[\widetilde{a},\widetilde{b}]=0$ and we put
$\pi_{2n+1} =\chi((g_{2n+1}, g_{2n+2}))$. Now  the properties of
the sequence imply that
$$
[\widetilde{g_1}, \widetilde{ g_2} \widetilde{g_3}] =
[\widetilde{g_{1}},\widetilde{ g_2}]\widetilde{g_3}=
\widetilde{g_3}\widetilde{g_2}\widetilde{g_1}(\pi_1\widetilde{b_1}-\widetilde{1}),
$$
and
\[
\begin{split}
[\widetilde{g_1}, & \widetilde{ g_2} \widetilde{g_3}, \widetilde{
g_4} \widetilde{g_5}, \ldots, \widetilde{g_{2n}}
\widetilde{g_{2n+1}}]=\\
& = \widetilde{g_{2n+1}}\widetilde{ g_{2n}}\cdots
\widetilde{g_{2}} \widetilde{g_{1}}
(\pi_{2n-1}\widetilde{b_{2n-1}}-\widetilde{1})\cdots
(\pi_3\widetilde{b_3}-\widetilde{1})
(\pi_1\widetilde{b_1}-\widetilde{1}).
\end{split}
\]
Clearly, this   is nonzero for each $n$, contradicting that
$F^\lambda [G]$ is upper (lower) Lie nilpotent. Therefore $G'$ is
a finite $p$-group, as claimed.

Finally, let $A'$ be of order $p^t$. Our assertion is valid for
$t=0$ and assume its truth for $t-1$. Lemma \ref{L:3} says that
$b^{p^t}$ belongs to the center of $G$ for any $b\in G$. It
follows that any conjugacy class of $G$, which belongs to $A'$,
has  $p$-power order or it is central. This yields that $A'$ has
central subgroup $L=\gp{c}$ of order $p$ and by Lemma \ref{L:2}
there exists $\gamma\in F$ such that $\widetilde{c}\gamma$ has
order $p$.  Then $F^\lambda [G]$ can be  realized in a second way
as a twisted group algebra $F^\tau [G]$ with the following
diagonal change of the basis
$$
\overline{g}= \begin{cases}
\widetilde{c}^i\gamma^i,                 &  \text{if }\quad  g=c^i;\\
\widetilde{g} &  \text{if }\quad  g\in G\setminus L
\end{cases}
$$
with   twisting  function  $\tau$. Since $\overline{c^i}$ is
central, Lemma \ref{L:5} asserts that
$$
F^\tau [G]/{\bf I}(L)\cong F^\tau [G/L].
$$
Of course, $F^\tau [G/L]$ is upper (lower) Lie nilpotent, so we
can  apply induction to obtain that $G'$ is a finite $p$-group.

\bigskip

{\it The second proof  for $char(F)>2$}.\quad
For any $a,b,c \in G$ the Lie commutator
$[\widetilde{a},\widetilde{b}, \widetilde{c}]$ belongs to the
nilpotent ideal $F^\lambda_\sigma [G]^{[3]}$ and let $x$ be a
non-zero fixed element  from  the annihilator of $F^\lambda_\sigma
[G]^{[3]}$. Clearly,
\begin{equation}\label{E:13}
\widetilde{a}\widetilde{b}
\widetilde{c}x-\widetilde{b}\widetilde{a} \widetilde{c}x=
\widetilde{c}\widetilde{a}\widetilde{b}x-
\widetilde{c}\widetilde{b} \widetilde{a}x,
\end{equation}
and without loss of generality we  can assume that $1 \in
Supp(x)$. Assume that  $Supp(x)=\{x_1=1, x_2,\ldots,x_l\}$. It is
convenient to distinguish  the  following  cases:

{\bf 1.}\quad $Supp(\widetilde{c}\widetilde{a}\widetilde{b}x)\cap
Supp(\widetilde{c}\widetilde{b} \widetilde{a}x)$ is not empty.
Then  $cabx_i=cbax_j$ for suitable $i$ and $j$, so the commutator
$(a,b)$ can be written as $x_jx_i^{-1}$.

{\bf 2.}\quad $Supp(\widetilde{c}\widetilde{a}\widetilde{b}x)\cap
Supp(\widetilde{c}\widetilde{b} \widetilde{a}x)$\quad  is an empty
set. The length of the right side of (\ref{E:13}) is $2l$ and thus
the length of the left one is also $2l$. This means that
$$
Supp(\widetilde{a}\widetilde{b} \widetilde{c}x)\cap
Supp(\widetilde{b}\widetilde{a} \widetilde{c}x)=\oslash.
$$
Now assume that\quad  $Supp(\widetilde{a}\widetilde{b}
\widetilde{c}x)\cap
Supp(\widetilde{c}\widetilde{a}\widetilde{b}x)$ and
$Supp(\widetilde{a}\widetilde{b} \widetilde{c}x)\cap
Supp(\widetilde{b}\widetilde{a} \widetilde{c}x)$\quad  are not
empty sets. There exist  $x_i,x_j,x_t,x_r\in Supp(x)$  such that
$abcx_i=cabx_j$,\quad $abcx_t=cbax_r$ and the commutator $(a,b)$
coincides with an element of the form $x_{i_1}^{\pm 1}x_{i_2}^{\pm
1}x_{i_3}^{\pm 1}x_{i_4}^{\pm 1}.$
 Similar statement is valid if
 $$
 Supp(\widetilde{b}\widetilde{a}
\widetilde{c}x)\cap
Supp(\widetilde{c}\widetilde{a}\widetilde{b}x)\quad \text{ and
}\quad  Supp(\widetilde{b}\widetilde{a} \widetilde{c}x)\cap
Supp(\widetilde{b}\widetilde{a} \widetilde{c}x)
$$
are not empty sets. It remains to consider one of the following
two subcases:

{\bf 2.1} \quad $\widetilde{a}\widetilde{b} \widetilde{c}x=
 \widetilde{c}\widetilde{a}\widetilde{b}x$ and
$\widetilde{b}\widetilde{a} \widetilde{c}x=
 \widetilde{c}\widetilde{b}\widetilde{a}x$. Then the commutators
$(ab,c)$ and $(ba,c)$ can be written as $x_jx_i^{-1}$.

{\bf 2.2}\quad  $\widetilde{a}\widetilde{b} \widetilde{c}x=
 -\widetilde{c}\widetilde{b}\widetilde{a}x$ and
$ \widetilde{c}\widetilde{a}\widetilde{b}x=-
\widetilde{b}\widetilde{a} \widetilde{c}x$. This  yields  that
$$
abcx_i=cbax_j\quad  \text{and}\quad  (ab,c)=x_jx_i^{-1}.
$$

Now, if $p\ne 2$, we  put $c=b^{-1}$ and then
$(ab,c)=(a,b)^{-1}$. The foregoing immediately implies that there
are only finitely many group commutators of the form $(a,b)$ with
$a,b\in G$; each commutator $(a,b)$ is a $p$-element which has
only  a finite number of conjugates. Then, as well known,  $G'$ is a finite $p$-group.

\bigskip

\begin{theorem} \label{T:1}
Let $F^\lambda_{\sigma} [G]$ be a crossed product of a group $G$
and the field $F$ of characteristic $0$ or $p$. Then
\begin{itemize}
\item[(1)] Any upper (lower) Lie nilpotent crossed product
$F^\lambda_{\sigma} [G]$ is a twisted group algebra. \item[(2)]
The twisted group algebra $F^\lambda[G]$ is lower (upper) Lie
nilpotent if and only if one of the following condition holds:
\begin{itemize}
\item[(2.i)] $F^\lambda [G]$ is a commutative algebra (i.e $G$ is
abelian and the twisting  function is symmetric). \item[(2.ii)]
$char(F)=p$, $G$ is a nilpotent group with commutator subgroup of
 $p$-power order and  the untwisted $p$-elements of $G$ form
a subgroup.  Moreover,  for any $a,b\in G$ the group   commutator
$(a,b)$ is an untwisted  $p$-element and
\newline
$\bigl(\lambda(a,b)^{-1}\lambda(b,b^{-1}ab)\lambda(a,(a,b)
\lambda((b,a),(a,b))^{-1}\bigr)^{-p^m}=\mu((a,b)), $
\newline
where  $p^m$ is the order of $(a,b)$.
\end{itemize}
\end{itemize}
\end{theorem}
{\bf Proof.}\quad  Let  $F^\lambda_{\sigma} [G]$ be a lower
(upper) Lie nilpotent crossed product of  characteristic 0 or $p$
and let $H=ker \sigma.$ The twisted group algebra $F^\lambda [H]$
is a subring of the crossed product $F^\lambda_{\sigma} [G]$ and,
by \cite{bb}, for every non-zero ideal $I$ of $F^\lambda_{\sigma}
[G]$ we have
\begin{equation}\label{E:12}
F^\lambda [H]\cap I \ne 0.
\end{equation}
Recall   that if  $char(F)=p$ and $\Delta(G)$ has no element of
order $p$, then Theorem 3.5 from \cite{mpass} states that
$F^\lambda_{\sigma} [G]$ is a semiprime ring. For $char F=0$
according to Corollary 6 from \cite{pas} the algebra
$F^\lambda[H]$ is semiprime and (\ref{E:12}) ensures   that
$F^\lambda_{\sigma} [G]$ is also semiprime. It follows that  $
F^\lambda_\sigma [G]^{[3]}=0$ and  (\ref{E:7})  implies that
$([a,b]r)^2=0$ for all $a,b,r\in F^\lambda_{\sigma} [G]$. Clearly,
$[a,b]F^\lambda_{\sigma} [G]$ is a nilpotent ideal, but this is
possible only if $F^\lambda_\sigma [G]^{[2]}=0$ and then
$F^\lambda_\sigma [G]$ is a commutative twisted group algebra, as required.

Finally assume that $p$ divides  the order of some
element of $\Delta(G)$ and  $F^\lambda_\sigma [G]$ is a
noncommutative crossed product. Then
 the Lie commutator
$
[\widetilde{a},\widetilde{1}\alpha]=\widetilde{a}(\alpha-\alpha^{\sigma(a)})
$ belongs to  the nil ideal $F^\lambda_\sigma [G]^{[2]}$ for every
$\alpha\in F$ and the element $\alpha-\alpha^{\sigma(a)}$ of the
field $F$  is zero, because it is nilpotent. Thus  $\sigma$ is
trivial and so the crossed product $F^\lambda_\sigma [G]$ is a
twisted group algebra.

Let  $F^\lambda [G]$ be a twisted group algebra
such that $char(K)=p$ and $F^\lambda [G]^{[p^t]}=0$.  By Lemma \ref{L:2}, the
commutator subgroup  $G'$ is a $p$-group and  Lemma \ref{L:7}
forces that $G'$ is finite. Furthermore,  Lemma \ref{L:3} says
that $b^{p^t}$ belongs to the center of $G$ for every $b\in G$. So
the quotient $G/C$ of $G$ by the center $C$ is a $p$-group of
finite exponent with finite commutator subgroup. Clearly, the
orders of those conjugacy classes of the group $G/C$, which are
contained in the finite normal $p$-subgroup $(G/C)'$,  are
$p$-powers. Hence $(G/C)'$ has a nontrivial central subgroup $L$
and thus $G/C$ is a nilpotent group. Remark that Lemma \ref{L:2}
confirms the remaining  statements.

The converse  statement was proved in \cite{bk}, so the proof is
complete.

\bigskip

Using the notation of untwisting,  Theorem \ref{T:1} can be
formulated as

\begin{corollary}\label{C:1}
Let $F^\lambda [G]$ be a twisted group algebra of a group $G$ and
a field $F$ of characteristic $0$ or $p>0$. The algebra
$F^\lambda[G]$ is lower (upper) Lie nilpotent if and only if one
of the following conditions holds:
\begin{itemize}
\item[(i)]  $F^\lambda [G]$ is  commutative; \item[(ii)]
$char(F)=p$, $G$ is a nilpotent group such that  $G'$ is a finite
$p$-group and  $F^\lambda[G]$ is  stably untwisted.
\end{itemize}
\end{corollary}

{\bf Proof.} Let  $F^\lambda [G]$ be a noncommutative lower
(upper) Lie nilpotent algebra.  By  Theorem \ref{T:1},
$char(F)=p$, $G$ is a nilpotent group and  $G'$  has $p$-power
order. As we remarked before, the  nilpotent elements of
$F^\lambda [G]$ form an ideal $N$ and $F^\lambda [G]/N$ is
commutative. By Lemma \ref{L:2} for any $g$ of $G'$ we  can choose
$\gamma_g\in F$ such that the order of $\widetilde{g}\gamma_g$
coincides with the order of $g$. Then
$\widetilde{g}\gamma_g-\widetilde{1}$ is nilpotent and belongs to
the  radical $J(F^\lambda [G'])$ of $F^\lambda [G']$. So
$J(F^\lambda [G'])$ is nilpotent of  codimension 1 and from
$$
J(F^\lambda [G'])\cdot F^\lambda [G]= F^\lambda [G]\cdot
J(F^\lambda [G'])
$$
 follows that $J(F^\lambda [G'])\cdot F^\lambda [G]$ is a
nilpotent ideal. For all $a,b\in G$ we have
\begin{equation}\label{E:21}
[\widetilde{a},\widetilde{b}]=\widetilde{a}^{-1}\widetilde{b}^{-1}
((\widetilde{a},\widetilde{b})-\widetilde{1})=\widetilde{a}^{-1}\widetilde{b}^{-1}(\widetilde{(a,b)}\chi((a,b))-
\widetilde{1}),
\end{equation}
and if $p^m$ is the order of the group commutator $(a,b)$ then by
Lemma 4 we obtain
$$
\bigl(\widetilde{(a,b)}\chi((a,b))-\widetilde{1}\big)^{p^m}=\mu((a,b))\chi((a,b))^{p^m}-\widetilde{1}=0.
$$
Consequently  $[\widetilde{a},\widetilde{b}]\in J(F^\lambda
[G'])\cdot F^\lambda [G]$, which  implies that
$$
 F^\lambda [G]/J(F^\lambda
[G'])\cdot F^\lambda [G]
$$
is a commutative algebra. Now Lemma \ref{L:1} states that
$F^\lambda [G]$ is  stably untwisted.

Conversely, if $F^\lambda [G]$ is stably untwisted, then there
exists an  extension $K$ of the field $F$ such that $K^\lambda
[G]$ is isomorphic to the ordinary group algebra $KG$ via a
diagonal change of basis and the result for $KG$ has  been already
known. Since $F^\lambda [G]$ is a subalgebra of $K[G]$, it is
lower (upper) Lie nilpotent.

\bigskip

\section{$(n,m)$-Engel crossed products}

Let $R$ be an associative ring and let $n,m$ be  fixed positive
integers. If  $$ [a,\underbrace{b^m,b^m,\ldots,b^m}_n]=0
$$
for all elements $a,b\in R$, then  $R$ is called  {\bf
$(n,m)$-Engel.}

Clearly, an  $(n,m)$-Engel ring  satisfies the polynomial identity
$$
[x,\underbrace{y^m,y^m,\ldots,y^m}_n].
$$
Let $p^t$ be the smallest positive integer such that $n\leq p^t$
and let $m=p^lr$ with $(p,r)=1$. If $F\gp{x,y}$ is a
noncommutative polynomial ring with indeterminates  $x$ and $y$
over a field $F$ of characteristic $p$, then
\begin{equation}\label{E:14}
\begin{split}
[x,\underbrace{y^m,y^m,\ldots,y^m}_{p^t}]&=[x,\underbrace{y^{p^lr},y^{p^lr},\ldots,y^{p^lr}}_{p^t}]\\
&=[x,y^{p^{l+t}r}]= [x,\underbrace{y^r,y^r,\ldots,y^r}_{p^{l+t}}].
\end{split}
\end{equation}
Therefore, if $R$ is an  $(n,m)$-Engel ring of $char(R)=p>0$, then it is $(p^{l+t},r)$-Engel ring, too.

\begin{theorem}\label{T:2}
Let $F^\lambda_{\sigma} [G]$ be a crossed product of a group $G$
and a  field $F$ of positive characteristic $p$.
\begin{itemize}
\item[(1)] Any  $(n,m)$-Engel crossed product $F^\lambda_{\sigma}
[G]$ is a twisted group algebra.

\item[(2)] If $F^\lambda [G]$ is  an $(n,m)$-Engel twisted group
algebra, then either $F^\lambda [G]$ is commutative, or the
following conditions hold:
\begin{itemize}
\item[(2.i)] $G$ has a normal subgroup $B$ of a finite  index such
that  commutator subgroup  $B'$  has   $p$-power order, the
$p$-Sylow subgroup $P/B$ of $G/B$ is a normal subgroup, $G/P$ is a
finite abelian group of an  exponent that divides $m$ and $P$ is a
nilpotent subgroup. \item[(2.ii)] the untwisted $p$-elements of
$G$ form a subgroup and  for all $a,b\in G$  the  commutator
$(a,b)$ is an untwisted $p$-element such that
\newline
$\bigl(\lambda(a,b)^{-1}\lambda(b,b^{-1}ab)\lambda(a,(a,b)
\lambda((b,a),(a,b))^{-1}\bigr)^{-p^m}=\mu((a,b))$,
\newline
where  $p^m$ is the order of $(a,b)$. Moreover, $F^\lambda[B]$ is
stably untwisted and $|G:B||B'|$ is bounded by a fixed function of
$n$ and $m$.
\end{itemize}
\end{itemize}
\end{theorem}
{\bf Proof.} Let $F^\lambda_\sigma[G]$ be an  $(n,m)$-Engel
crossed product  of characteristic $p>0$. Then  we can  apply
Theorem 3 of Kezlan \cite{kez} which states  that
$F^\lambda_\sigma [G]^{[2]}$ is a nil ideal. Hence the nilpotent
elements of $F^\lambda_\sigma [G]$ form an ideal  $N$ and
$F^\lambda_\sigma [G]/N$ is commutative. Clearly,\quad
$[\widetilde{a},\widetilde{1}\cdot
\alpha]=\widetilde{a}(\alpha-\alpha^{\sigma(a)})$\quad  belongs to
the nil ideal $F^\lambda_\sigma [G]^{[2]}$ for every $\alpha\in
F$. It follows that the element $\alpha-\alpha^\sigma(a)$ of the
field $F$  is zero, because it is nilpotent. So   $\sigma$ is
trivial and  hence  $F^\lambda_\sigma [G]=F^\lambda [G]$ is a
twisted group algebra.

First assume that $G$ has no element of order $p$. By Corollary 2
from \cite{passm},  the nil ideal $F^\lambda[G]^{[2]}$ is zero. So
the twisted group algebra $F^\lambda [G]$ is commutative.

Now let $G$ be a noncommutative group with  a $p$-element. Without
loss of generality, by (14), we   can  assume that $(m,p)=1$.
Lemma \ref{L:2} ensures that every group commutator and their
products are $p$-elements. Therefore $G'$ is a $p$-group. Choose
the smallest positive integer $p^t$ with $n\leq p^t$ and let $a,
b\in G$. The   Lie commutator formula ensures that
$$
[\widetilde{a},\widetilde{b}^m,p^t]=[\widetilde{a},
\underbrace{\widetilde{b}^m,\widetilde{b}^m,\ldots,
\widetilde{b}^m}_{p^t}]=[\widetilde{a},\widetilde{b}^{mp^t}]=0.
$$
This  yields
$\widetilde{a}\widetilde{b}^{mp^t}=\widetilde{b}^{mp^t}\widetilde{a}$.
Hence $b^{mp^t}$ is  central for any $b\in G$. By  Lemma
\ref{L:11},\quad  $G$ has a normal subgroup $A$ of finite index
such that commutator subgroup $A'$ has  $p$-power order. Our aim
is to prove that if $B$ is the subgroup generated by $A$ and the
center of $G$, then $F^{\lambda}[B]$ is stably untwisted. Of
course,  $B'$ is a finite $p$-group, so  the ideal $I=J(F^\lambda
[B'])\cdot F^\lambda [B]$ is nilpotent, because the nilpotent
elements of $F^\lambda [G]$ form an ideal. Indeed, if $p^r$ is the
order of the commutator $(c,d)$ of $B$, then Lemma \ref{L:2}
ensures that $((\widetilde{c},\widetilde{d})-1)^{p^r}=0$\quad and
$$
(\widetilde{c},\widetilde{d})-1\in J(F^\lambda [B'])\cdot
F^\lambda [B].
$$
Then   (\ref{E:21}) asserts  that\quad
$[\widetilde{c},\widetilde{d}]\in  J(F^\lambda [B'])\cdot
F^\lambda [B]$,\quad  and \quad $ F^\lambda [B]/J(F^\lambda
[B'])\cdot F^\lambda [G]$\quad  is a commutative algebra. Now, by
Lemma \ref{L:1}, $F^\lambda_{\sigma} [B]$ is  stably untwisted, as
required.

Now $b^{mp^t}$  is central for any $b\in G$ and it belongs to $B$
and the $p$-Sylow subgroup $P/B$ of $G/B$ is normal, because the
commutator subgroup of $G$ is a $p$-group. Then, by the
Schur-Zassenhaus's theorem (6.2.1 theorem of \cite{gor}), $G/B$ is
a semidirect product of  the finite $p$-Sylow subgroup $P/B$ and a
finite abelian group $M/B$ whose  exponent divides $m$. Since
$B^{mp^t}$ is a central subgroup,  there remains  to show that
$P/B^{mp^t}$ is nilpotent. To prove this we proceed by induction
on the order of the commutator subgroup $L$ of $B/B^{mp^t}$. First
suppose that $L'=\gp{1}$. The  abelian group $B/B^{mp^t}$ is a
direct product of the  $p$-subgroup $D$ and a subgroup $H$ whose
exponent divides $m$. The finite $p$-group $P/B$ acts by
conjugation  on these subgroups. Of course for every $d\in
P/B^{mp^t}$ and $h\in H$ the $p'$-element $d^{-1}hd$ coincides
with $h(h,d)$ and the commutator $(h,d)$ is a $p$-element.
Therefore, the action of $P/B$ on $H$ is trivial and $P/B$ acts as
a finite $p$-group on $D$.

Define the subgroups\quad  $D_1=(P/B^{mp^t},D)$\; and \;
$D_j=(P/B^{mp^t},D_{j-1})$\quad  for $j>1$. By  Lemma V.4.1
\cite{sh}, this sequence of  subgroups is such that $D_l=\gp{1}$
for some $l$. But the finite $p$-group $P/B$ is nilpotent, so a
suitable term $K_s$ of the lower central series of the group
$K=P/B^{mp^t}$ is contained in $B/B^{mp^t}$ and
$$
K_{s+1}=(K,K_s)\subseteq (K,B/B^{mp^t})= (K,D\times H)
\subseteq (K,D)=D_1,
$$
because  $H$ is central in $P/B^{mp^t}.$ Similarly, we conclude
that $K_{s+i}\subseteq D_i$ and the group $P/B^{mp^t}$ must be
nilpotent, because  $D_l=1$.

Finally,  assume that  $L'\not=\gp{1}$. The orders of those
conjugacy classes of the group $P/B^{mp^t}$, which are contained
in the finite normal $p$-subgroup $L'$, are $p$-powers. Hence $L'$
has a nontrivial central  subgroup, and by induction on the order
of $L'$, we see that $P/B^{mp^t}$ is nilpotent. Consequence, $P$
is nilpotent as asserted.

\begin{corollary}\label{C:2}
The twisted group algebra  $F^\lambda [G]$ of positive
characteristic $p$ is  $n$-Engel if and only if either $F^\lambda
[G]$ is commutative, or the following conditions hold:
\begin{itemize}
\item[(i)]   $G$ is a nilpotent group with  a normal subgroup $B$
of a finite  $p$-power index, $B'$ is a finite $p$-group and
$F^\lambda[B]$ is stably untwisted; \item[(ii)]  the untwisted
$p$-elements of $G$ form a subgroup, the commutator $(a,b)$ is an
untwisted $p$-element for all $a,b\in G$ and
\begin{equation}\label{E:15}
\bigl(\lambda(a,b)^{-1}\lambda(b,b^{-1}ab)\lambda(a,(a,b)
\lambda((b,a),(a,b))^{-1}\bigr)^{-p^m}=\mu((a,b)),
\end{equation}
 where  $p^m$ is the order of $(a,b)$.
\end{itemize}
\end{corollary}

{\bf Proof.} There remains to prove only the sufficiency of these
conditions, and as before,  we can assume  that $F$ is an
algebraically closed field of characteristic $p$.

Let $G$ be a nilpotent group with a normal abelian  subgroup $B$
of a finite $p$-power index. Then for any $a,b$ of $B$ the twist
of $\widetilde{(a,b)}$ is equal to 1. So (\ref{E:15}) asserts that
$F^\lambda [B]$ is a commutative algebra. Now we try to adapt the
method of the  proof of  Theorem V.6.1 from \cite{sh}. For this we
need additional information about the nilpotent groups $G$ with a
normal abelian subgroup $B$ of index $p^s$ and about the twisted
group algebras. Certainly, $(B,G^{p^s})=1$ and Lemma V.6.2 from
\cite{sh} implies that $(B, G)^{p^m}=1$ for some $m\geq s$. It
follows that
$$
(G^{p^{s+m}},G) \subseteq (B^{p^m},G)=\gp{1}.
$$
Hence, if $t=s+m$, then $g^{p^t}$ is central in $G$ and
(\ref{E:15}) confirms that $[\widetilde{g^{p^t}},\widetilde{h}]=0$
for any $h\in G$. So $\widetilde{g}^{p^t}$ is central. Now for
every $y=\sum_{g\in G}c_g\widetilde{g}$ of $F^\lambda [G]$ we have
that \quad $\sum_{g\in G} c_g^{p^t} \widetilde{g}^{p^t}$\quad  is
central in $F^\lambda [G]$ and
$$\textstyle
y^{p^t}= \sum_{g\in G} c_g^{p^t} \widetilde{g}^{p^t} + y_1,
$$
for suitable $ y_1\in [F^\lambda [G],F^\lambda [G]]$.   For all
$a,b\in G$ we have
\begin{equation}\label{E:22}
[\widetilde{a},\widetilde{b}]=\widetilde{a}^{-1}\widetilde{b}^{-1}
((\widetilde{a},\widetilde{b})-\widetilde{1})=\widetilde{a}^{-1}\widetilde{b}^{-1}
(\widetilde{(a,b)}\chi((a,b))-\widetilde{1})
\end{equation}
and if $p^m$ is the order of $(a,b)$ then by (\ref{E:15}),
$$
\bigl(\widetilde{(a,b)}\chi((a,b))-\widetilde{1}\big)^{p^m}=\mu((a,b))\chi((a,b))^{p^m}-\widetilde{1}=0.
$$
Consequently  $[F^\lambda [G],F^\lambda [G]] \subseteq
F^\lambda[G]\cdot J(F^\lambda [G'])$. Since the subgroup $D=(G,
B)$ is  normal in $G$ and $B/D$ is central in  $G/D$ of index\quad
$|G/D:B/D|=p^s$,\quad  it follows  that $G'/D$ is a finite group
of order $p^l$. By Theorem 1.6 \cite{passr}, we have
$$
(F^\lambda[G]\cdot J(F^\lambda [G'])^{p^l}\subseteq
J(F^\lambda[D])\cdot F^\lambda [G].
$$
It follows that\quad  $\textstyle y_1^{p^l}= \sum_ {i=1}^{p^s} z_i
\widetilde{t_i}$,\quad   where $t_1=1, t_2,\ldots,t_{p^s}$ is a
transversal of $B$ in $G$. Furthermore,
$$
z_i\in F^\lambda[B]\cap J(F^\lambda[D])\qquad (i=1, 2,\ldots, p^s)
$$
are nilpotent elements with  nilpotency index at most $p^m$ and
each of these elements commute, because $B$ is abelian. The inner
automorphisms of $\widetilde{G}$ induces a finite group of
automorphisms $T$ on $\widetilde{B}$. The action of $T$ on $z_1,
z_2,\ldots,z_{p^s}$ produces only finitely many images and denote
by $L$ the subring of  the commutative
 algebra $F^\lambda[B]$ generated by these images. Of course  $L$ is nilpotent, its
nilpotency index is at most\quad  $p^r=p^{m+s+1}|T|$\quad  which
does not depend on $L$. Clearly, $y_1^{p^{l+r}}=0$, and by the
foregoing,
$$
\textstyle y^{p^{t+l+r}}= \sum_{g\in G} c_g^{p^{t+l+r}}
\widetilde{g}^{p^{t+l+r}}
$$
is central in $FG$. By the   Lie commutator identity we obtain
that
$$[x,y,p^{t+l+r}]=[x,y^{p^{t+l+r}}]=0$$ and  $FG$ is Lie
$p^{t+l+r}$-Engel.

Finally, let $B'$ be of order $p^t$. Our assertion is valid for
$t=0$, assume its truth for $t-1$. The normal subgroup $B'$ of a
nilpotent group  contains a central  cyclic subgroup
$L=\gp{c\mid\,c^p=1}$. Now we take $F^\lambda G$ and make the
following  change of basis:
$$
\overline{g}= \begin{cases}
\widetilde{g},          quad        &  \text{if }\quad  g\in G\setminus L;\\
\widetilde{c^i}\sqrt{{\mu(c)}^{-i}} \quad &  \text{if }\quad  g=c^i\in
L.
\end{cases}
$$
Then $  \rho (c^i,c^j)= 1$ and $\overline{c}^p=\widetilde{1}$.
 By Lemma \ref{L:5}, $F^\lambda[G]$  can be realized in a second way as a twisted group ring
with new basis $\{\overline{\overline{g}}\}$ and a twisting
function $\tau (g,h)$ which satisfies $\tau(c^k,g)=1$.
Clearly,\quad  $\tau(c^k,g)=\tau(g,c^k)=1$.\quad   By
(\ref{E:24}), we have that
\begin{equation*}
F^\lambda G/{\bf I}(L)\cong F^\tau [G/L]
\end{equation*}
and by (\ref{E:4}),\quad   $ \tau(g_iH,g_jH)=\lambda(g_i,g_j)$.\quad Hence
the twisting  function $\tau(a,b)$ satisfies the required
conditions, and by the inductive hypothesis,  $F^\mu [G/L]$ is
$p^m$-Engel for some $m$. It follows that \quad
$[x,y,p^m]=[x,y^{p^m}]\in {\bf I}(L)$\quad for all $x,y\in
F^\lambda G$, so  \quad $[x,y^{p^m}]=(\overline{c}-1)z$\quad for
some $z\in F^\lambda G$. Since $\overline{c}$ is central, we have
$$\textstyle
[x,y^{p^m},p]= [x,\underbrace{y^{p^m},y^{p^m},\ldots
,y^{p^m}}_p]\in (\overline{c}-1)^p F^\lambda G =0.
$$
This implies  that $[x,y,p^{m+1}]=0$ and the proof is complete.

\bibliographystyle{abbrv}
\bibliography{Bovdi_Grishkov_final}

\begin{thebibliography}{10}

\bibitem{bb}
A.~A. Bovdi.
\newblock Cross products of semigroups and rings.
\newblock {\em Sibirsk. Mat. Z.}, 4:481--499, 1963.

\bibitem{bk}
A.~A. Bovdi and K.~K. Kolikov.
\newblock Lie nilpotency and ideals of crossed products.
\newblock {\em Serdica}, 15(4):275--286 (1990), 1989.

\bibitem{gor}
D.~Gorenstein.
\newblock {\em Finite groups}.
\newblock Chelsea Publishing Co., New York, second edition, 1980.

\bibitem{gupt}
N.~Gupta and F.~Levin.
\newblock On the {L}ie ideals of a ring.
\newblock {\em J. Algebra}, 81(1):225--231, 1983.

\bibitem{kez}
T.~P. Kezlan.
\newblock Some rigs with nil commutator ideals.
\newblock {\em Michigan Math. J.}, 12:105--111, 1965.

\bibitem{liu}
C.-H. Liu.
\newblock On units of twisted group algebras.
\newblock {\em J. Algebra}, 250(1):271--282, 2002.

\bibitem{mpass}
S.~Montgomery and D.~S. Passman.
\newblock Crossed products over prime rings.
\newblock {\em Israel J. Math.}, 31(3-4):224--256, 1978.

\bibitem{neu}
B.~H. Neumann.
\newblock Groups with finite classes of conjugate elements.
\newblock {\em Proc. London Math. Soc. (3)}, 1:178--187, 1951.

\bibitem{pps}
I.~B.~S. Passi, D.~S. Passman, and S.~K. Sehgal.
\newblock Lie solvable group rings.
\newblock {\em Canad. J. Math.}, 25:748--757, 1973.

\bibitem{passr}
D.~S. Passman.
\newblock Radicals of twisted group rings.
\newblock {\em Proc. London Math. Soc. (3)}, 20:409--437, 1970.

\bibitem{pass}
D.~S. Passman.
\newblock Group rings satisfying a polynomial identity. {II}.
\newblock {\em Pacific J. Math.}, 39:425--438, 1971.

\bibitem{passm}
D.~S. Passman.
\newblock Trace methods in twisted group algebras.
\newblock {\em Proc. Amer. Math. Soc.}, 129(4):943--946, 2001.

\bibitem{pas}
D.~S. Passman.
\newblock Twisted group algebras satisfying a generalized polynomial identity.
\newblock {\em Comm. Algebra}, 29(9):3683--3710, 2001.
\newblock Special issue dedicated to Alexei Ivanovich Kostrikin.

\bibitem{sh}
S.~K. Sehgal.
\newblock {\em Topics in group rings}, volume~50 of {\em Monographs and
  Textbooks in Pure and Applied Math.}
\newblock Marcel Dekker Inc., New York, 1978.

\end{thebibliography}

\end{document}